\documentclass[12pt]{article}
\usepackage{amssymb,amsmath,latexsym,mathrsfs,eufrak,graphics,graphicx,longtable}
\newtheorem{theorem}{Theorem}

\newtheorem{proposition}{Proposition}
\usepackage[left=1cm, top=1cm, bottom=2.5cm, right=1cm]{geometry}
\begin{document}
\def\Pr{{\bf pr}}
\def\Ad{{\bf Ad }}
\def\di{\displaystyle }
\title{PRELIMINARY GROUP CLASSIFICATION AND SOME EXACT SOLUTIONS OF $2-$HESSIAN EQUATION}
\author{Mahdieh Yourdkhany \thanks{Department of Mathematics Karaj branch of Islamic Azad University, Karaj, Iran. e-mail:~m.yourdkhany@kiau.ac.ir} \and Mehdi Nadjafikhah \thanks{School of Mathematics, Iran University of Science and Technology, Narmak, Tehran, Iran. e-mail:~m\_nadjafikhah@iust.ac.ir} \and Megerdich toomanian \thanks{(Corresponding author)Department of Mathematics Karaj branch of Islamic Azad University, Karaj, Iran. email:~megerdich.toomanian@kiau.ac.ir}}
\maketitle
\abstract{We study the class  of $3$-dimensional nonlinear $2-$hessian equations $u_{xx}u_{yy}+u_{xx}u_{yy}+u_{yy}u_{zz}-u_{xy}^2-u_{yz}^2-u_{xz}^2-f(x,y,z)=0$, where $f$ is an arbitrary smooth function of the variables $(x,y,z)$. We perform preliminary group classification on $2-$hessian equation. In fact, we find additional equivalence transformation on the space $(x,y,z,u,f)$, with the aid of \textit{N. Bila's method}, then we take their projections on the space $(x,y,z,f)$, so we prove an optimal system of one-dimensional Lie subalgebras of this equation is generated by $\langle A^1,\cdots, A^{12}\rangle$, which introduced in theorem  \eqref{opti}, ultimately, A number of new interesting nonlinear invariant models are obtained which have non-trivial invariance algebras. The result of these works is a wide class of equations which summarized in table. So at the end of this work, some exact solutions of $2-$hessian equation are presented. The paper is one of the few applications of an algebraic approach to the group classification using Lie method.}

$\textrm{Keywords:}$ Hessian equation, Optimal system, Preliminary group classification

$\textrm{AMS Classification 2010:}$ 53C10, 53C12, 53A55, 35A30, 76M60, 58J70

\section{Introduction}
Nowadays it is generally accepted that a huge number of real processes arising in physics, biology, chemistry, etc. can be described by nonlinear PDEs. And the most powerful methods for costruction of exact solutions for a wide ranges of nonlinear PDEs are symmetry-based methods, and these methods originated from the Lie method, so the basic part of the theory, is infinitesimal method of Sophus Lie~\cite{Lie}, that is connection between continuous transformation groups and algebras of their infinitesimal generators. This method leads to techniques in the group-invariant solutions and conservation laws of differential  equations~\cite{O1,Ib3,Ov1}.

In fact the method that we proposed, the method of preliminary group classification, is a conclusion of Lie infinitesimal  method, and is defined  and related to the theory of group classification of differential equations. this method is proposed in~\cite{Akhatov} and is developed for deferential equation in~\cite{Ib2,Bihlo}.
 
The main idea of preliminary group classification is based on extension of the kernel of admitted Lie groups that 
are obtained by the transformations from the corresponding equivalence Lie group. The problem of finding inequivalent cases of such extension of symmetry can reduce to the classification of inequivalent subgroups of the equivalence Lie group(In particular, if a Lie group is finite-parameter, then one can use an optimal systems of its subgroups). we use equivalence transformations and the theory of classification of finite-dimensional Lie algebras. In this paper we study point symmetry and equivalence classification of $\mathrm{HESI}$ equation, leading to a number of new interesting nonlinear invariant models associated to non-trivial invariance algebras.

A complete list of these models are given for a finite-dimensional equivalence algebra derived for $\mathrm{HESI}$ equation. To obtain these goals we perform algorithms that is explained in references ~\cite{Ov1,Bila,Ib4,Cherniha}, and we use  similar works in ~\cite{Karn,Ib2,mah1,Dr1,Dr2,Ib1,Song}.

\textbf{For the local solution:} The existence of  $\textit{C}^{\infty}$ local solutions of $\mathrm{HESI}$ equation In $\mathbb{R}^3$ is studied in~\cite{Tian}, and the solution is in the following form: 
\begin{equation*}
u(x,y,z)=\frac{1}{2}(\tau_1x^2+\tau_2y^2+\tau_3z^2)+\varepsilon^5\omega(\varepsilon^{-2}(x,y,z)), 
\end{equation*}
where $\varepsilon$ and $\tau_i$  are arbitrary constants, and $\omega$ is a given smooth function. So at the end of this work, with the symmetry group of the equation, this solutions  transforms to another solution of $\mathrm{HESI}$ equation.

\textbf{About $\mathrm{HESI}$ equation:} Based on refrences~\cite{fros,Tian}, $k$-Hessian equations are a family of PDEs in $n$-dimensional space equations that can be written as ${\cal S}_k[u]=f$, where $1\leqslant k \leqslant n$, ${\cal S}_k[u]=\sigma_k(\lambda({\cal D}^2u))$, and $\lambda({\cal D}^2u)=(\lambda_1,\cdots,\lambda_n)$, are the eigenvalues of the \textit{Hessian matrix} ${\cal D}^2u$ $\big((\partial_i \partial_ju)_{1\leqslant i,j \leqslant n}\big)$ and  $\sigma_k(\lambda)=\sum_{i_1<\cdots<i_k}\lambda_{i_1}\cdots\lambda_{i_k}$, is a $k$th elementry symmetric polynomial.

The  $k$-Hessian equations include the Laplace equation, when $k=1$, And the Monge-Ampere equation, when $k=n$. 

Here we study $2-$hessian equation in three dimensions, and f is an arbitrary function of $x,y,z$ ($\mathrm{HESI}$ equation) :
\begin{equation}
{\cal S}_2[u] := u_{xx}u_{yy}+u_{xx}u_{yy}+u_{yy}u_{zz}-u_{xy}^2-u_{yz}^2-u_{xz}^2,
\end{equation}
this equation is a fully nonlinear elliptic partial diferential equation, that is related to intrinsic curvature for three-dimensional manifolds.

In fact, the $2-$hessian equation is unfamiliar outside Riemannian geometry and elliptic regularity theory, that is closely related to the scalar curvature operator, which provides an intrinstic curvature for a three-dimensional manifold. Geometric PDEs have been used widely in image analysis ~\cite{sapiro}. In particular, the Monge-Ampere equation in the context of optimal transportation has been used in three dimensional volume-based image registration ~\cite{haker}. 

The $2-$hessian operator  also appears in conformal mapping problems. Conformal surface mapping has been used for two-dimmensional image registration ~\cite{angenent,Gu}, but does not generalize directly to three dimensions. Quasi-conformal maps have been used in three dimensions ~\cite{wang,zeng}. However, these methods are still being developed.
\section{Principal Lie Algebra}
The symmetry approach to the classification of admissible partial differential equations depends heavily on a useful way of describing transformation groups that keep invariant the form of a given partial differential equation. This is done via the well-known infinitesimal method developed by Sophus Lie~\cite{O1,O2,Ov1}.Given a partial differential equation, the problem of finding its maximal (in some sense) Lie invariance group reduces to solving \textit{determining equations} that is an over-determined system of linear partial differential equation.
We consider the $2-$hessian equation as the form: 
 \begin{equation}
\mathrm{HESI}:\; {\cal S}_2[u]=f(x,y,z),\label{Hesi}
 \end{equation}
where $u=u(x,y,z)$ is dependent variable and $x,y,z$ are independent variables and $f$ is arbitrary function. Considering the total space $ E=X\times U$ with local coordinate $(x,y,z,u)$ which $x,y,z\in X$ and $u\in U$. The solution space of equation (\ref{Hesi}) is a subvariety $S_{\Delta}\subset J^2(\mathbb{R}^3,\mathbb{R})$ of the second order of jet bundle of $3-$dimensional submanifolds of $E$. The $1-$parameter Lie group of infinitesimal transformations on $E$ is as follows:
\begin{equation}
\begin{aligned}
\tilde{x}&=x+t\xi(x,y,z,u)+O(t^2), & \tilde{y}&=y+t\zeta(x,y,z,u)+O(t^2), \\ 
\tilde{z}&=z+t\eta(x,y,z,u)+O(t^2), & \tilde{u}&=u+t\phi(x,y,z,u)+O(t^2), 
\end{aligned} 
\end{equation}
where $t$ is the group parameters and $\xi,\zeta,\phi$ and $\eta$ are the infinitesimals of the transformations for the independent and dependent variables, resp. So the corresponding infinitesimal generators have the following form generally 
\begin{equation}
V=\xi(x,y,z,u)\partial_x +\zeta\partial_y +\eta\partial_z +\phi\partial_u.
\end{equation}
So based on Theorem 2.31 of~\cite{O2}, $V$ is a invariant point transformation if $\Pr^{(2)}V[\mathrm{HESI}]=0$. Where $\Pr^{(2)}V$ is the second order prolongation of the vector field $V$, that means:
\begin{equation}
 \Pr^{(2)}V = V+{\phi}^x(x,y,z,u^{(2)})\,{\partial}_{u_x}+
 \cdots+{\phi}^{xz}(x,y,z,u^{(2)})\,{\partial}_{u_{xz}},
\end{equation}
 in which $u^{(2)}=(u,u_x,u_y,u_z,u_{xx},u_{xy},u_{xz},u_{yy},u_{yz},u_{zz})$ and
 \begin{align}
 {\phi}^J(x,y,z,u^{(2)})=\textbf{D}_J\bigg(\phi-\sum_{i=1}^{3}{\xi}^i\frac{\partial u}{\partial x_i}\bigg)+\sum_{i=1}^{3}{\xi}^i\frac{\partial {u_J}}{\partial x_i}
 \end{align}
 where $({\xi}^1,{\xi}^2,{\xi}^3)=(\xi,\zeta,\eta)$ and $(x_1,x_2,x_3)=(x,y,z)$, further $J=(j_1,\cdots,j_k)$ is a $k$-th order multi-index, and $j_i$s adopt $x$ or $y$ or $z$, for each $1\leqslant i \leqslant k$, then $\textbf{D}_J$ denotes the total derivatives for the multi-index $J$, and the $J$-th total derivative is as 
 $\textbf{D}_J=\textbf{D}_{j_1}\cdots\textbf{D}_{j_k}$, where
\begin{equation*}
\textbf{D}_i={\partial}_{x_i}+\sum_{J}\frac{\partial {u_J}}{\partial x_i}\,{\partial}_{u_J}, \qquad (x_1,x_2,x_3)=(x,y,z).
\end{equation*} 
So $\Pr^{(2)}V$ acts on Eq.\eqref{Hesi} and with replacing $u_{yy}$ with equivalent expression of $\mathrm{HESI}$ equation we have the following system as determining equation: 
\begin{eqnarray}
\begin{aligned}
&{\phi}_{xx}= {\phi}_{xy}= {\phi}_{xz}={\phi}_{xu}= {\phi}_{yy}= {\phi}_{yz}={\phi}_{yu}= {\phi}_{zz}= {\phi}_{zu}=0,\\
&{\phi}_{uu}= {\xi}_{u}= {\xi}_{yy}={\xi}_{yz}= {\xi}_{zz}= {\zeta}_{u}={\zeta}_{zz}= {\eta}_{u}= {\eta}_{zz}=0,\\
&{\xi}_{x}={\eta}_{z}, \;\; {\zeta}_{x}=-{\xi}_{y}, \;\; {\zeta}_{y}={\eta}_{z},\;\;{\eta}_{x}=-{\xi}_{z}, \;\; {\eta}_{y}=-{\zeta}_{z}, \\
&f_x{\xi}+f_y{\zeta}+f_z{\eta}+2f(2{\eta}_z-{\phi}_u)=0.
\end{aligned}
\end{eqnarray}
where $f$ is arbitrary function.with solving above relations we have:
\begin{eqnarray} \label{con}
\begin{aligned}
&\xi=c_6x+c_7y+c_8z+c_9, \qquad \zeta=c_{10}z+c_6y-c_7x+c_{11},\\
&\eta=-c_{10}y+c_6z-c_8x+c_{12}, \qquad \phi=c_1x+c_2u+c_3y+c_4z+c_5, \\
&f_x{\xi}+f_y{\zeta}+f_z{\eta}+2f(2{\eta}_z-f{\phi}_u)=0.
\end{aligned}
\end{eqnarray}
which $c_i,i=1,\cdots,12$ are arbitrary costants.

So if $f(x,y,z)=0$ the last equation of relations \eqref{con} will be removed and based on the first four equations in relations \eqref{con}, we have $12$-dimensional symmetry group, but if $f(x,y,z)\neq0$ we substitute the first four equations of \eqref{con} in the last one, and obtain the following \textit{condition}:
\begin{eqnarray} \label{con2}
\begin{aligned}
(c_6x+c_7y&+c_8z+c_9)f_x+( c_{10}z+c_6y-c_7x+c_{11})f_y \\
&+(-c_{10}y+c_6z-c_8x+c_{12})f_z+(-2c_2+4c_6)f=0.
\end{aligned} 
\end{eqnarray}
There isn't $c_1,c_3,c_4,c_5$ in condition \eqref{con2}, so these coeficients are free,that means equation \eqref{Hesi} have $4-$dimentional symmetry group at minimum. So we conclude the following theorem from above relations:
\begin{theorem}:
The $\mathrm{HESI}$ equation(Eq.\eqref{Hesi}) admits symmetry group of dimension 4 to 12, for different choises of given function $f(x,y,z)$. These equations have the common following vectors as infinitesimal generators:
\begin{align}\label{principle}
V_1&=\partial_u , & V_2&=x\partial_u , & V_3&=y\partial_u , & V_4&=z\partial_u.
\end{align}
\end{theorem}

Then the Lie algebra $\mathfrak{g}$ generated with the vectors \eqref{principle} is called the principal Lie algebra for Eq.\eqref{Hesi}. Now we want to specify the coefficient $f$ such that Eq.\eqref{Hesi} admits an extension of the principal algebra $\mathfrak{g}$. therefore, we do not solve the determining equation, instead we obtain a partial group classification of Eq.\eqref{Hesi} via so-called method of  \textit{preliminary group classification}.

This method was suggested in ~\cite{Akhatov} and applied when an equivalence group is generated by a finite-dimensional Lie algebra $\mathfrak{g}_{\mathscr{E}}$. The essential part of the method is the classification of all nonsimilar subalgebras of $\mathfrak{g}_{\mathscr{E}}$. Actually the classification is based on finite-dimentional equivalence algebra $\mathfrak{g}_{\mathscr{E}}$.
\section{Equivalence Transformations}

with a nondegenerate change of the variables $x,y,z$ an equation of the form $\mathrm{HESI}$ equation convert to an equation of the same form, but with different $f(x,y,z)$. The set of all equivalence transformatioms forms an equivalence group $E$. We shal find a subgroup $E_c$ of it with infinitesimal method.

We suppose an operator of the group $E_c$ is in the form:
 \begin{equation}
 Y=\xi(x,y,z,u)\partial_x +\zeta\partial_y +\eta\partial_z +\phi\partial_u +\psi(x,y,z,u,f)\partial_f.\label{Y}
 \end{equation}
So from the invariance conditions of Eq.\eqref{Hesi} written as the following system:
\begin{align}
{\cal S}_2[u]=f(x,y,z),\qquad f_u=0.
\end{align}

Note that $f$ and $u$ are considered as differential variables; $u$ on the space $(x,y,z)$ and $f$ on the space $(x,y,z,u)$. The coordinates $\xi,\zeta,\eta,\phi$ of operator \eqref{Y} are funtions of $x,y,z,u$, while the coordinate $\psi$ is function of $x,y,z,u,f$. as usual way we should solve the following system that obtained of the invariance conditions:
\begin{align}
\Pr^{(2)}Y( {\cal S}_2[u]=f(x,y,z)),\qquad \Pr^{(2)}Y(f_u)=0.
\end{align}
Where $\Pr^{(2)}Y$ is the second order prolongation of the vector field $Y$.
 
 \textbf{But}, to obtain the operator $Y$ of the group $E_c$ we use of \textit{N. Bila's method} in ref.~\cite{Bila}. The base of our procedure is theorem (1) of ~\cite{Bila}, then this theorem and it's results can be summarized as the following three-steps procedure:
 
 \textbf{step 1:} Find the determining equations of the extended classical symmetries related to the Eq.\eqref{Hesi}. For the meaning of the extended classical symmetries, a vector 
 \begin{equation}
 V=\xi(x,y,z,u,f)\partial_x +\zeta\partial_y +\eta\partial_z +\phi\partial_u +\psi\partial_f.\label{v}
 \end{equation}
 is said the extended classical symmetry operator assosiated with $\mathrm{HESI}$ Equation and the determining equations of the extended classical symmetries related to the $\mathrm{HESI}$ Equation is the following equation:
 \begin{equation}
 \Pr^{(2)}V[\mathrm{HESI}]=0,\label{sys}
 \end{equation}
 where $\xi$, $\zeta$, $\eta$, $\phi$ and $\psi$ are functions of $x$, $y$, $z$, $u$ and $f$, and $\Pr^{(2)}V$ is
 \begin{align} 
 V+\sum_{J}{\phi}^J(x,y,z,u^{(2)},f^{(2)})\,{\partial}_{u_J}+\sum_{J}{\psi}^J(x,y,z,u^{(2)},f^{(2)})\,{\partial}_{f_J}
 \end{align}
where $u^{(2)}=(u, u_x,\cdots,u_{zz})$ and $f^{(2)}=(f, f_x,\cdots, f_{zz})$: 
 and the coefficients obtain from:
 \begin{align} 
 {\phi}^J(x,y,z,u^{(2)},f^{(2)})=\textbf{D}_J\bigg(\phi-\sum_{i=1}^{3}{\xi}^i\frac{\partial u}{\partial x_i}\bigg)+\sum_{i=1}^{3}{\xi}^i\frac{\partial {u_J}}{\partial x_i}\nonumber\\
 {\psi}^J(x,y,z,u^{(2)},f^{(2)})=\textbf{D}_J\bigg(\psi-\sum_{i=1}^{3}{\xi}^i\frac{\partial f}{\partial x_i}\bigg)+\sum_{i=1}^{3}{\xi}^i\frac{\partial {f_J}}{\partial x_i}
 \end{align}
 where $({\xi}^1,{\xi}^2,{\xi}^3)=(\xi,\zeta,\eta)$ and $(x_1,x_2,x_3)=(x,y,z)$, further $J=(j_1,\cdots,j_k)$ is a $k$-th order multi-index, and $j_i$s adopt $x$, $y$ or $z$, for $1\leqslant i \leqslant k$, then $\textbf{D}_J$ denotes the total derivatives for the multi-index $J$, and the $J$-th total derivative is as 
$\textbf{D}_J =\textbf{D}_{j_1}\cdots\textbf{D}_{j_k}$, that total derivative operator with respect to $i$ is as following
 \begin{equation*}
\textbf{D}_i=\partial_{x_i}+\sum_{J}\frac{\partial {u_J}}{\partial x_i}\,{\partial}_{u_J}+\sum_{J}\frac{\partial {u_f}}{\partial x_i}\,{\partial}_{f_J},
 \end{equation*}
\textbf{Note:} $\Pr^{(2)}V$ is determined by taking into account that $u$ and $f$ are both dependent variables, exactly as one would proceed in finding the classical Lie symmetries for a system without arbitrary functions.

So with solving equation \eqref{sys} we gain:
\begin{equation}\label{sys2}
\begin{aligned}
& \xi=c_6x+c_9y+c_7z+c_8, & & \zeta=c_{10}z+c_6y-c_9x+c_{11}, \\
& \psi=2f(-2c_6+c_3), & & \eta=-c_{10}y+c_6z-c_7x+c_{12}, \\
& \phi=c_1x+c_3u+c_2y+c_5z+c_4,
\end{aligned}
\end{equation}
which $c_i,i=1,\cdots,12$ are arbitrary costants.

\textbf{step 2:} Augment the system of step 1 with the following conditions:
\begin{equation}
\frac{\partial {\xi}}{\partial u}=0,\qquad \frac{\partial {\zeta}}{\partial u}=0,\qquad\frac{\partial {\eta}}{\partial u}=0,\qquad \frac{\partial {\psi}}{\partial u}=0.
 \end{equation}
As we seen in relations \eqref{sys2} above conditions are satiesfied.

\textbf{step 3:} Augment the system of steps 1 and 2 with the following conditions:
\begin{equation}
\frac{\partial {\xi}}{\partial f}=0,\qquad \frac{\partial {\zeta}}{\partial f}=0,\qquad\frac{\partial {\eta}}{\partial f}=0,\qquad \frac{\partial {\psi}}{\partial f}=0.
 \end{equation}
As we seen in relations \eqref{sys2} above conditions are satisfied too.

Ultimately, The class of equations \eqref{Hesi} has a finite continuous group of equivalence transformations generated by the following infinitesimal operators:
\begin{equation} \label{vectors2}
\begin{aligned}
Y_1 &=\partial_x , & Y_2 &=\partial_y , & Y_3 &=\partial_z ,\\
Y_4 &=\partial_u , & Y_5 &=x\partial_u , & Y_6 &=y\partial_u ,\\
Y_7 &=z\partial_u , & Y_8 &=z\partial_x -x\partial_z , & Y_9 &=y\partial_x -x\partial_y ,\\
Y_{10} &=z\partial_y -y\partial_z , & Y_{11} &=u\partial_u +2f\partial_f, & Y_{12} &=x\partial_x +y\partial_y +z\partial_z -4f\partial_f.
\end{aligned}
\end{equation}
Moreover, in the group of equivalence transformations are included also discrete transformations, i.e., reflections
$(x,y,z,u,f)\mapsto-(x,y,z,u,f)$.
\section{Sketch of the method of preliminary group classification}
In many applications of group analysis, most of extensions of the principal Lie algebra admitted by an equation are obtained from the equivalence algebra $\mathfrak{g}_{\mathscr{E}}$. We call these extension ${\mathscr{E}}$-\textit{extension of the principal Lie algebra}. The classification of all nonequivalent equations admitting ${\mathscr{E}}$-extension of the principal Lie algebra is called \textit{a preliminary qroup classification}. What we obtain also is not necessarily the largest equivalence group but, it can be any subgroup of the qroup of all equivalence transformations.

The application of this method is effective and simple when it is based on a finite-dimensional equivalence algebra $\mathfrak{g}_{\mathscr{E}}$. So we take finite dimensional algebra $\mathfrak{g}_{12}$ spanned on the basis \eqref{vectors2} and use it for preliminary group classification.

The function $f$ of Eq.\eqref{Hesi} depends on the variables $x,y,z$, so we don't construct any prolongations of operators \eqref{Y}. But we take projections on the space $(x,y,z,f)$.

The nonzero projections of \eqref{vectors2} are:
\begin{align}
Z_1 &=\Pr(Y_1)=\partial_x , & Z_2 &=\Pr(Y_2)=\partial_y ,\nonumber\\
Z_3 &=\Pr(Y_3)=\partial_z , & Z_4 &=\Pr(Y_8)=z\partial_x -x\partial_z , \nonumber\\
Z_5 &=\Pr(Y_9)=y\partial_x -x\partial_y , & Z_6 &=\Pr(Y_{10})z\partial_y -y\partial_z ,\nonumber\\
Z_7 &=\Pr(Y_{11})=2f\partial_f, & Z_8 &=\Pr(Y_{12})=x\partial_x +y\partial_y +z\partial_z -4f\partial_f.\label{vectors3}
\end{align}
It's clear that there aren't the minimal infinitesimal generators \eqref{principle}, among above vectors.

The Lie algebra generated with the vectors in \eqref{vectors3} is denoted by $\mathfrak{g}_{8}$. 

The essence of the preliminary method is based on the following two proposition:

\begin{proposition}:\label{proposition1}
Let $\mathfrak{g}_{m}$ be a $m-$dimensional subalgebra of $\mathfrak{g}_{8}$. Suppose $Z^{(i)}$, $i=1,\cdots,m$ be a basis of $\mathfrak{g}_{m}$ and $Y^{(i)}$ is the elements of the algebra $\mathfrak{g}_{12}$, such that $Z^{(i)}=\Pr(Y^{(i)})$, that means, if
\begin{equation}\label{1}
 Z^{(i)}=\sum_{\alpha=1}^{8}e^\alpha_iZ_{\alpha},
 \end{equation}
 then with respect to \eqref{vectors2} and \eqref{vectors3}:
 \begin{align}\label{2}
 Y^{(i)}=e^{1}_iY_{1}+e^{2}_iY_{2}+e^{3}_iY_{3}+e^{4}_iY_{8}+e^{5}_iY_{9}+e^{6}_iY_{10}+
 e^{7}_iY_{11}+e^{8}_iY_{12}.
\end{align} 
If function $f=f(x,y,z)$ be invariant with respect to the algebra $\mathfrak{g}_{m}$, then the $\mathrm{HESI}$ equation admits the operators
 \begin{equation}\label{3}
 X^{(i)}=\mbox{projection of}\;\;Y^{(i)} \mbox{on}\;\; (x,y,z,u).
\end{equation}
\end{proposition}
\begin{proposition}:\label{proposition2}
Let equations 
\begin{eqnarray}
 {\cal S}_2[u]=f(x,y,z), \label{eq1} \\
{\cal S}_2[u]=f^{\prime}(x,y,z),\label{eq2}
\end{eqnarray}
 be constructed according to proposition \eqref{proposition1} with subalgebras $\mathfrak{g}_{m}$ and $\mathfrak{g}_{m^{\prime}}$, respectively. If $\mathfrak{g}_{m}$ and $\mathfrak{g}_{m^{\prime}}$ are similar subalgebras in $\mathfrak{g}_{12}$ then equations \eqref{eq1} and \eqref{eq2} are equivalent with respect to the equivalence group $G_{12}$ generated by $\mathfrak{g}_{12}$. 
\end{proposition}

According to above propositions, continuation of the preliminary group classification of Eq.\eqref{Hesi} with respect to the finite-dimensional algebra $\mathfrak{g}_{12}$, is reduced to the algebraic problem of constructing of nonsimilar subalgebras of $\mathfrak{g}_{8}$, or optimal systems of subalgebras. 

\textbf{note:} In this paper we just solve the problem of preliminary group classification with respect to one-dimensional subalgebras. 
\section{Adjoint group for algebra $\mathfrak{g}_{8}$ }
We determine a list or optimal system, of conjuacy inequivalent subalgebras with the property that any other subalgebra is equivalent to a unique member of the list under some element of the adjoint representation, i.e. $\bar{\mathfrak{h}}=\Ad(g)\mathfrak{h}$ for some g of a considered Lie group, see~\cite{O1,O2,Ov1}. 

The adjoint action is given by the Lie series 
 \begin{equation}
 \Ad(\exp(\varepsilon Y_{i}))Y_{j}=Y_{j}-\varepsilon [Y_{i},Y_{j}]+\frac{\varepsilon^2}{2}[Y_{i},[Y_{i},Y_{j}]]-\cdots,
 \end{equation}
The commutator and adjoint representations of $\mathfrak{g}_{8}$ are listed in tables \ref{tab1} and \ref{tab2}. 

\begin{table}[ht]
\caption{Commutators table for $\mathfrak{g}_{8}$: $[Z_{i},Z_{j}]$} \label{tab1}
\begin{center}
\begin{tabular}{c|cccccccc}
 & $Z_{1}$ & $Z_{2}$ & $Z_{3}$ & $Z_{4}$ & $Z_{5}$ & $Z_{6}$ & $Z_{7}$ & $Z_{8}$\\
\hline
$Z_{1}$ & 0 & 0 & 0 & $-Z_{3}$ & $-Z_{2}$ & 0 & 0 & $Z_{1}$\\
$Z_{2}$ & 0 & 0 & 0 & 0 & $Z_{1}$ & $-Z_{3}$ & 0 & $Z_{2}$\\
$Z_{3}$ & 0 & 0 & 0 & $Z_{1}$ & 0 & $Z_{2}$ & 0& $Z_{3}$\\
$Z_{4}$ & $Z_{3}$ & 0 & $-Z_{1}$ & 0 & $-Z_{6}$ & $Z_{5}$ & 0 & 0\\
 $Z_{5}$ & $Z_{2}$ & $-Z_{1}$ & 0 & $Z_{6}$ & 0 & $-Z_{4}$ & 0 & 0\\
$Z_{6}$ & 0 & $Z_{3}$ & $-Z_{2}$ & $-Z_{5}$ & $Z_{4}$ & 0 & 0 & 0\\
$Z_{7}$ & 0 & 0 & 0 & 0 & 0 & 0 & 0 & 0\\
 $Z_{8}$ & $-Z_{1}$ & $-Z_{2}$ & $-Z_{3}$ & 0 & 0 & 0 & 0 & 0\\
\end{tabular}
\end{center}
\end{table}
\begin{table}[ht]
\begin{center}
\caption{Adjoint table for $\mathfrak{g}_{8}$: $\Ad(\exp(\varepsilon_i Y_{i}))Y_{j}$} \label{tab2}
\scalebox{0.8}{
\begin{tabular}{c|cccccccc}
  & $Z_{1}$ & $Z_{2}$ & $Z_{3}$ & $Z_{4}$ & $Z_{5}$ & $Z_{6}$ & $Z_{7}$ & $Z_{8}$\\
\hline
$Z_{1}$ & $Z_{1}$ & $Z_{2}$ & $Z_{3}$ & $\varepsilon_1Z_{3}+Z_{4}$ & $\varepsilon_1Z_{2}+Z_{5}$ & $Z_{6}$ & $Z_{7}$ & $-\varepsilon_1Z_{1}+Z_{8}$\\
\hline
$Z_{2}$ & $Z_{1}$ & $Z_{2}$ & $Z_{3}$ & $Z_{4}$ & $-\varepsilon_2Z_{1}+Z_{5}$ & $\varepsilon_2Z_{3} +Z_{6}$ & $Z_{7}$ & $-\varepsilon_2Z_{2}+Z_{8}$\\
\hline
$Z_{3}$ & $Z_{1}$ & $Z_{2}$ & $Z_{3}$ & $-\varepsilon_3Z_{1}+Z_{4}$ & $Z_{5}$ & $-\varepsilon_3Z_{2} +Z_{6}$ & $Z_{7}$ & $-\varepsilon_3Z_{3}+Z_{8}$\\
\hline
$Z_{4}$ & $\di \cos(\varepsilon_4)Z_{1} \atop \di-\sin(\varepsilon_4)Z_{3}$ & $Z_{2}$ & $\di \sin(\varepsilon_4)Z_{1} \atop \di  +\cos(\varepsilon_4)Z_{3}$ & $Z_{4}$ & $\di\cos(\varepsilon_4)Z_{5} \atop \di+\sin(\varepsilon_4)Z_{6}$ & $\di-\sin(\varepsilon_4)Z_{5} \atop \di +\cos(\varepsilon_4)Z_{6}$ & $Z_{7}$ & $Z_{8}$ \\
\hline
 $Z_{5}$ & $\di \cos(\varepsilon_5)Z_{1} \atop \di -\sin(\varepsilon_5)Z_{2}$ & $\di \sin(\varepsilon_5)Z_{1} \atop \di +\cos(\varepsilon_5)Z_{2}$ & $Z_{3}$ & $\di \cos(\varepsilon_5)Z_{4} \atop \di-\sin(\varepsilon_5)Z_{6}$ & $Z_{5}$ & $\di \sin(\varepsilon_5)Z_{4} \atop \di +\cos(\varepsilon_5)Z_{6}$ & $Z_{7}$ & $Z_{8}$ \\
 \hline
$Z_{6}$ & $Z_{1}$ & $\di \cos(\varepsilon_6)Z_{2} \atop \di -\sin(\varepsilon_6)Z_{3}$ & $\di \sin(\varepsilon_6)Z_{2} \atop \di  +\cos(\varepsilon_6)Z_{3}$ & $\di \cos(\varepsilon_6)Z_{4} \atop \di +\sin(\varepsilon_6)Z_{5}$ & $\di -\sin(\varepsilon_6)Z_{4} \atop \di  +\cos(\varepsilon_6)Z_{5}$ & $Z_{6}$ & $Z_{7}$ & $Z_{8}$ \\
\hline
$Z_{7}$ & $Z_{1}$ & $Z_{2}$ & $Z_{3}$ & $Z_{4}$ & $Z_{5}$ & $Z_{6}$ & $Z_{7}$ & $Z_{8}$\\
\hline
$ Z_{8}$ & $e^{\varepsilon_8}Z_{1}$ & $e^{\varepsilon_8}Z_{2}$ & $e^{\varepsilon_8}Z_{3}$ & $Z_{4}$ & $Z_{5}$ & $Z_{6}$ & $Z_{7}$ & $Z_{8}$
\end{tabular}}
\end{center}
\end{table}
\section{Construction of the optimal system of one-dimensional subalgebras of $\mathfrak{g}_{8}$}
\begin{theorem}:\label{opti}
An optimal system of one-dimensional Lie algebras of $\mathfrak{g}_{8}$ in $\mathrm{HESI}$ equation are as follows:
\begin{equation}\label{op}
\begin{aligned}
& A^1=Z_{7}, & A^2&=\pm Z_{1}+Z_{7},\\
& A^3=\gamma_1Z_{6}+Z_{7}, & A^4&=\pm Z_{1}+\gamma_2Z_{6}+Z_{7},\\
& A^5=\alpha_1Z_{4}+Z_{7}, & A^6&=\pm Z_{2}+\alpha_2Z_{4}+Z_{7},\\
&A^7=\alpha_3Z_{4}+\gamma_3Z_{6}+Z_{7}, & A^8&=\pm Z_{1}+\alpha_4Z_{4}+\gamma_4Z_{6}+Z_{7},\\
& A^9=\alpha_5Z_{4}+\beta_1Z_{5}+Z_{7}, & A^{10}&=\pm Z_{3}+\alpha_6Z_{4}+\beta_2Z_{5}+Z_{7},\\
& A^{11}=\alpha_7Z_{4}+\beta_3Z_{5}+\gamma_5Z_{7}+Z_{8}, & A^{12}&=\pm Z_{2}+\alpha_8Z_{4}+\beta_4Z_{5}+\gamma_6Z_{7}+Z_{8},
\end{aligned}
\end{equation}
Where $\alpha_i$, $i=1,...8$ and $\beta_i$, $i=1,\cdots,4$ and $\gamma_i$, $i=1,\cdots,6$ are arbitrary constants.
\end{theorem}
\medskip \noindent {\it Proof:} 
We will start with $Z=\sum_{i=1}^{8}a_iZ_{i}$, suppose Z is a nonzero vector field of $\mathfrak{g}_{8}$, we want simplify as many of the coefficients $a_i$, $i=1,\cdots,8$ as possible through proper Adjoint applications on Z. We proceed this simplifications through following cases:

\textbf{note:} The coefficients $a_7$ and $a_8$ don't change at all.

\textbf{Case 1:} At first, we assume that $a_8=0$, so with scalling on Z, we can suppose that $a_7=1$, then we have $Z=\sum_{i=1}^{6}a_iZ_{i}+Z_7$. therefore for different values of $a_5=0$, when it is either zero or nonzero, we have cases 1.1 and 1.2.

\textbf{Case 1.1:} If $a_8=a_5=0$, so $Z=a_1Z_{1}+a_2Z_{2}+a_3Z_{3}+a_4Z_{4}+a_6Z_{6}+Z_{7}$. Then for different values of $a_4$, when it is either zero or nonzero, we have cases 1.1.a and 1.1.b.

\textbf{Case 1.1.a:} If $a_8=a_5=a_4=0$, so $Z=a_1Z_{1}+a_2Z_{2}+a_3Z_{3}+a_6Z_{6}+Z_{7}$. Then for different values of $a_6$, when it is either zero or nonzero, we have cases 1.1.a1 and 1.1.a2.

\textbf{Case 1.1.a1:} If $a_8=a_5=a_4=a_6=0$, so $Z=a_1Z_{1}+a_2Z_{2}+a_3Z_{3}+Z_{7}$. Then for different values of $a_3$, when it is either zero or nonzero, the coefficient can be vanished; when $a_3\neq 0$, with effecting $\Ad(\exp(\cot^{-1}(a_1/{a_3})Z_4))$ on Z. Then we have $Z=a_1Z_{1}+a_2Z_{2}+Z_{7}$.

 Now if $a_2=0$ or $a_2\neq0$; by effecting $\Ad(\exp(\cot^{-1}(a_1/{a_2})Z_5))$ on Z, we can make the coefficient of $Z_2$ vanished.Then we have $Z=a_1Z_{1}+Z_{7}$.
 
 So if $a_1=0$, then $Z=Z_7$, so we have $A^1$.

 And if $a_1\neq 0$, with $\Ad(\exp(\ln(\pm1/{a_1})Z_8))$ change the coefficient of $Z_1$ equal $\pm 1$, so $Z=\pm Z_1+Z_7$, therefore we have $A^2$.
 
 \textbf{Case 1.1.a2:} If $a_8=a_5=a_4=0$, but $a_6\neq0$, so $Z=a_1Z_{1}+a_2Z_{2}+a_3Z_{3}+a_6Z_{6}+Z_{7}$. Then for different values of $a_3$, when it is either zero or nonzero, the coefficient can be vanished; when $a_3\neq 0$, with applying $\Ad(\exp(-a_3/{a_6})Z_2)$ on Z. So we have $Z=a_1Z_{1}+a_2Z_{2}+a_6Z_{6}+Z_{7}$. 

 Similarly, the coefficient $a_2$ is either zero or we make it vanished with effecting $\Ad(\exp(a_2/{a_6})Z_3)$ on Z. Then we have $Z=a_1Z_{1}+a_6Z_{6}+Z_{7}$.
 
 So $a_1=0$ or $a_2\neq0$, if $a_1=0$, so $Z=a_6Z_{6}+Z_7$, and we have $A^3$.

 And if $a_1\neq 0$, with $\Ad(\exp(\ln(\pm1/{a_1})Z_8))$ change the coefficient of $Z_1$ equal $\pm 1$, so $Z=\pm Z_1++a_6Z_{6}+Z_7$, therefore we have $A^4$.
 
 \textbf{Case 1.1.b:} If $a_8=a_5=0$, but $a_4\neq0$, so $Z=a_1Z_{1}+a_2Z_{2}+a_3Z_{3}+a_4Z_{4}+a_6Z_{6}+Z_{7}$. Then for different values of $a_3$, when it is either zero or nonzero, the coefficient can be vanished; when $a_3\neq 0$, with effecting $\Ad(\exp(-a_1/{a_4})Z_1))$ on Z. Then we have $Z=a_1Z_{1}+a_2Z_{2}+a_4Z_{4}+a_6Z_{6}+Z_{7}$.

Therefore, for different values of $a_6$, when it is either zero or nonzero, we have cases 1.1.b1 and 1.1.b2.

\textbf{Case 1.1.b1:} If $a_8=a_5=a_3=a_6=0$, but $a_4\neq0$, so $Z=a_1Z_{1}+a_2Z_{2}+a_4Z_{4}+Z_{7}$. Then either $a_1$ is zero or nonzero, but
the coefficient can be vanished; when $a_1\neq 0$, with effecting $\Ad(\exp(a_1/{a_4})Z_3))$ on Z. Then we have $Z=a_2Z_{2}+a_4Z_{4}+Z_{7}$.

 Ultimately, $a_2=0$ or $a_2\neq0$; If $a_2=0$, so we have $Z=a_4Z_{4}+Z_{7}$.Then we have $A^5$.

 And if $a_2\neq 0$, with $\Ad(\exp(\ln(\pm1/{a_2})Z_8))$ change the coefficient of $Z_2$ equal $\pm 1$, so $Z=\pm Z_2+a_4Z_{4}+Z_7$, therefore we have $A^6$.
 
 \textbf{Case 1.1.b2:} If $a_8=a_5=a_3=0$, but $a_6\neq0$ and $a_4\neq0$, so $Z=a_1Z_{1}+a_2Z_{2}+a_4Z_{4}+a_6Z_{6}+Z_{7}$. Then for different values of $a_2$, when it is either zero or nonzero, the coefficient can be vanished; when $a_2\neq 0$, with applying $\Ad(\exp(a_2/{a_6})Z_3)$ on Z. So we have $Z=a_1Z_{1}+a_4Z_{4}+a_6Z_{6}+Z_{7}$. 

 Similarly, the coefficient $a_1$ is either zero or nonzero, if $a_1=0$, so $Z=a_4Z_{4}+a_6Z_{6}+Z_7$, and we have $A^7$.

 And if $a_1\neq 0$, with $\Ad(\exp(\ln(\pm1/{a_1})Z_8))$ change the coefficient of $Z_1$ equal $\pm 1$, so $Z=\pm Z_1+a_4Z_{4}++a_6Z_{6}+Z_7$, therefore we have $A^8$. 
 
 \textbf{Case 1.2:} If $a_8=0$ but $a_5\neq0$, so $Z=a_1Z_{1}+a_2Z_{2}+a_3Z_{3}+a_4Z_{4}+a_5Z_{5}+a_6Z_{6}+Z_{7}$. So we have different values of $a_2$, when it is either zero or nonzero, the coefficient can be vanished; when $a_2\neq 0$, with applying $\Ad(\exp(-a_2/a_5)Z_1)$ on Z. Then again we have two cases $a_1=0$ or $a_1\neq0$; If $a_1\neq0$, with applying $\Ad(\exp(a_1/a_5)Z_2)$ on Z, make it vanished. so we have $Z=a_3Z_{3}+a_4Z_{4}+a_5Z_{5}+a_6Z_{6}+Z_{7}$.
 
 Again we have two cases $a_6=0$ or $a_6\neq0$; If $a_6\neq0$, with applying $\Ad(\exp(\cot^{-1}(a_4/{a_6}))Z_5)$ on Z, make it vanished. so we have $Z=a_3Z_{3}+a_4Z_{4}+a_5Z_{5}+Z_{7}$. 
 
 Ultimately, If $a_3=0$, then $Z=a_4Z_{4}+a_5Z_{5}+Z_{7}$.Then we have $A^9$. 
 
 And if $a_3\neq0$, then with $\Ad(\exp(\ln(\pm1/{a_3})Z_8))$ change the coefficient of $Z_3$ equal $\pm 1$, so $Z=\pm Z_3+a_4Z_{4}+a_5Z_{5}+Z_{7}$.Then we have $A^{10}$.
 
 \textbf{Case 2:} If $a_8\neq0$, so with scalling on Z, we can suppose that $a_8=1$, then we have $Z=\sum_{i=1}^{7}a_iZ_{i}+Z_8$. So we have different values of $a_1$, when it is either zero or nonzero, the coefficient can be vanished; when $a_1\neq 0$, with applying $\Ad(\exp(a_1Z_1))$ on Z. So we reduce Z to $Z=a_2Z_{2}+a_3Z_{3}+a_4Z_{4}+a_5Z_{5}+a_6Z_{6}+a_7Z_{7}+Z_{8}$.
 
 Now if $a_3=0$ or $a_3\neq0$; by effecting $\Ad(\exp(\cot^{-1}(a_2/{a_3})Z_6)$ on Z, make the coefficient of $Z_3$ vanished. so we have $Z=a_2Z_{2}+a_4Z_{4}+a_5Z_{5}+a_6Z_{6}+a_7Z_{7}+Z_{8}$.
 
 Again we have two cases $a_6=0$ or $a_6\neq0$; If $a_6\neq0$, with applying $\Ad(\exp(-\cot^{-1}(a_5/{a_6}))Z_4)$ on Z, make it vanished. So we have $Z=a_2Z_{2}+a_4Z_{4}+a_5Z_{5}+a_7Z_{7}+Z_{8}$.
 
 Ultimately, If $a_2=0$, then $Z=a_4Z_{4}+a_5Z_{5}+a_7Z_{7}+Z_{8}$. Then we have $A^{11}$.
 
 And if $a_2\neq0$, then with $\Ad(\exp(\ln(\pm1/{a_2})Z_8))$ change the coefficient of $Z_2$ equal $\pm 1$, so $Z=\pm Z_2+a_4Z_{4}+a_5Z_{5}+a_7Z_{7}+Z_{8}$. Then we have $A^{12}$.
 
 There is not any more possible cases, and the proof is complete. \hfill\ $\Box$
\section{Equations admitting an extension by one of the principal Lie algebra}
Now based on propositions \eqref{proposition1} and \eqref{proposition2}, and with the optimal system \eqref{op}, we obtain all nonequivalent equations of the form equation \eqref{Hesi}, that admitting extension of principal Lie algebra $\mathfrak{g}$ by one operator $V_5$, that means every equation of the form equation \eqref{Hesi} admits symmetry group of dimension $4$ with infinitesimal generators \eqref{principle}, also together with a fifth operator $V_5$. for every case, when this extension occurs, we indicate the corresponding coefficients $f$ and additional operator $V_5$.
 
The algorithm will be clarified with these examples:
\paragraph{First example:} Consider the operator $A^3=\gamma_1Z_{6}+Z_{7}$, which $\gamma_1\neq 0$, from \eqref{op}, so \begin{equation}
\begin{aligned}
 A^3=\gamma_1z\partial_y -\gamma_1y\partial_z +2f\partial_f.
 \end{aligned}
 \end{equation}
 
 Invariants are found from the following equation (see ~\cite{O2}):
 \begin{equation}
 \begin{aligned}
 \frac{dy}{\gamma_1z}=-\frac{d z}{\gamma_1y}=\frac{d f}{2f},
 \end{aligned}
 \end{equation}
 and are the following functions:
 \begin{equation}
 \begin{aligned}
 I_1=x,\qquad I_2=y^2+z^2,\qquad I_3=f\exp(-\frac{2}{\gamma_1}{\tan}^{-1}(\frac{y}{z})).
 \end{aligned}
 \end{equation}
It follows
 \begin{equation}
 \begin{aligned}
f={\exp(\frac{2}{\gamma_1}{\tan}^{-1}(\frac{y}{z}))}{H(x,y^2+z^2)},
 \end{aligned}
 \end{equation}
 where $H$ is arbitrary function.
 
 By applying the formulas \eqref{1}, \eqref{2} and \eqref{3} on the operator $A^3$ we obtain the additional operator $V_5=\gamma_1z\partial_y -\gamma_1y\partial_z +u\partial_u$. Thus, the equation
 \begin{equation}
 \mathrm{HESI}: \; {\cal S}_2[u]={\exp(\frac{2}{\gamma_1}{\tan}^{-1}(\frac{y}{z}))}{H(x,y^2+z^2)},
 \end{equation}
admits the five-dimensional algebra ${\mathfrak{g}}_5$ , that is generated with the following vectors
 \begin{align}
V_1&=\partial_u , & V_2&=x\partial_u , &V_3&=y\partial_u , \nonumber\\ 
V_4&=z\partial_u , & V_5&=\gamma_1z\partial_y -\gamma_1y\partial_z +u\partial_u.
\end{align}
\paragraph{Second example:} Consider the operator $A^1=Z_{7}$ from \eqref{op}, so $A^1=2f\partial_f$. 
 
Invariants are the following functions:
 \begin{equation}\label{inv}
 \begin{aligned}
 I_1=x,\qquad I_2=y,\qquad I_3=z.
 \end{aligned}
 \end{equation}
 So, there are no invariant functions $f=f(x,y,z)$ because the necessary condition for existence of invariant solutions based on ref.~\cite{Ov1} (section 19.3) is not satisfied; that means invariants \eqref{inv} can't be solved with respect to $f$.
 
We continue calculations on some operators of \eqref{op}, and show results in table \ref{tab3}, that is the preliminary group classification of equation \eqref{Hesi}, which admit an extension ${\mathfrak{g}}_5$ of the principal Lie algebra ${\mathfrak{g}}$.

The results of classification: 
\begin{longtable}{c|cc}
\caption{The equation ${\cal E}:\;{\cal S}_2[u]=f$ has $V_5$ as its additional operator v.r.t $A^s$}  \label{tab3}\\
           & $f$ & $V_5$ \\ \hline 
$A^2$ & $ {\exp(\pm 2x)}{H(y,z)}$ &  $ \pm\partial_x +u\partial_u $\\
\hline 
${A^3}_{({\gamma_1}\neq 0)}$ & $\di\exp(\frac{2}{\gamma_1}{\tan}^{-1}(\frac{y}{z})) H(x,y^2+z^2)$ &  $ \gamma_1z\partial_y -\gamma_1y\partial_z +u\partial_u $\\
\hline 
${A^4}_{({\gamma_2}\neq 0)}$ &  $\di \exp(\frac{2}{\gamma_2}{\tan}^{-1}(\frac{y}{z})) \atop \di H(y^2+z^2,x\mp \frac{1}{\gamma_2}{\tan}^{-1}(\frac{y}{z}))$ &  $ \pm\partial_x +\gamma_2z\partial_y -\gamma_2y\partial_z +u\partial_u $\\
\hline 
${A^5}_{({\alpha_1}\neq 0)}$ &  $ \exp(\frac{2}{\alpha_1}{\tan}^{-1}(\frac{x}{z}))H(y,x^2+z^2)$ &  $  \alpha_1z\partial_x -\alpha_1x\partial_z +u\partial_u $\\
\hline 
${A^6}_{({\alpha_2}=0)}$ &  $ {\exp(\pm 2y)}{H(x,z)}$ & $ \pm\partial_y +u\partial_u $\\
\hline 
${A^6}_{({\alpha_2}\neq 0)}$ & $\di \exp(\frac{2}{\alpha_2}{\tan}^{-1}(\frac{x}{z})) \atop \di H(x^2+z^2,y\mp \frac{1}{\alpha_2}{\tan}^{-1}(\frac{x}{z}))$ &  $ \pm\partial_y +\alpha_2z\partial_x -\alpha_2x\partial_z +u\partial_u $\\
\hline 
${A^7}_{({\alpha_3}\neq 0,{\gamma_3}\neq 0)}$ &  $\di \exp(\frac{2}{\sqrt{{\alpha_3}^2+{\gamma_3}^2}}{\tan}^{-1}{(\frac{\alpha_3x+\gamma_3y}{z\sqrt{\alpha_3^2+\gamma_3^2}})}) \atop \di H(y-\frac{\gamma_3}{\alpha_3}x,(
(1-\frac{{\gamma_3}^2}{{\alpha_3}^2})x^2+\frac{2\gamma_3}{\alpha_3}xy+z^2))$ &  $\di \alpha_3z\partial_x +\gamma_3z\partial_y \atop \di  -(\alpha_3x+\gamma_3y)\partial_z +u\partial_u $\\
\hline 
${A^9}_{({\alpha_5}=0,{\beta_1}\neq 0)}$ &  $ \exp(\frac{2}{\beta_1}{\tan}^{-1}(\frac{x}{y}))H(z,x^2+y^2)$ &  $ \beta_1y\partial_x -\beta_1x\partial_y +u\partial_u$\\
\hline 
${A^9}_{(\alpha_5\neq 0,\beta_1\neq 0)}$ &  $\di \exp(\frac{-2}{\sqrt{\beta_1^2+\alpha_5^2}}\tan^{-1}(
\frac{\beta_1y+\alpha_5z}{x\sqrt{\alpha_5^2+\beta_1^2}})) \atop \di H(z-\frac{\alpha_5}{\beta_1}y, (x^2+(1-\frac{\alpha_5^2}{\beta_1^2})y^2+\frac{2\alpha_5}{\beta_1}yz))$ &  $\di (\alpha_5z+\beta_1y)\partial_x \atop \di  -\beta_1x\partial_y -\alpha_5x\partial_z +u\partial_u $\\
\hline 
 ${A^{10}}_{(\alpha_6=\beta_2=0)}$ &  $ {\exp(\pm 2z)}{H(x,y)}$ &  $ \pm\partial_z +u\partial_u $\\
 \hline 
${A^{10}}_{(\alpha_6\neq 0,\beta_2=0)}$ &  $\di \exp(\frac{2}{\beta_2}{\tan}^{-1}(\frac{x}{y})) \atop \di H(x^2+y^2,z\mp \frac{1}{\beta_2}{\tan}^{-1}(\frac{x}{y}))$ &  $ \pm\partial_z +\beta_2y\partial_x -\beta_2x\partial_y +u\partial_u $\\ 
\hline 
${A^{11}}_{(\alpha_7=\beta_3=\gamma_5=0)}$ &  $ x^{-4}H(\frac{y}{x},\frac{z}{x})$ &  $ x\partial_x +y\partial_y +z\partial_z $\\ 
\hline 
${A^{11}}_{(\alpha_7=\beta_3=0,\gamma_5\neq 0)}$ &  $ x^{2c-4}H(\frac{y}{x},\frac{z}{x})$ & $ x\partial_x +y\partial_y +z\partial_z $\\ 
\hline 
${A^{12}}_{(\alpha_8=\beta_4=\gamma_6=0)}$ &  $ x^{-4}H(\frac{y\pm 1}{x},\frac{z}{x})$ &  $ x\partial_x +(y\pm 1)\partial_y +z\partial_z $\\ 
\hline 
${A^{12}}_{(\alpha_8=\beta_4=0,\gamma_6\neq 0)}$ &  $ x^{2c-4}H(\frac{y\pm 1}{x},\frac{z}{x})$ &  $ x\partial_x +(y\pm 1)\partial_y +z\partial_z $ 
\end{longtable}
\section{Some Local Solutions }
Based on the following theorem of~\cite{Tian}, The $2-$hessian equation in $\mathbb{R}^3$,
\begin{equation}\label{2hesi}
{\cal S}_2[u]=f(x,y,z,u,{\cal D}u) \qquad  \mbox{on }  \Omega\subset\mathbb{R}^3,
\end{equation}
which $f\in\textit{C}^{\infty}(\Omega\times\mathbb{R}\times\mathbb{R}^3)$, ${\cal D}u=(\partial_1u,\cdots,\partial_nu)$, has $\textit{C}^{\infty}$ local solutions in $\mathbb{R}^3$, and the solution is in the following form: 
\begin{equation}\label{solution}
u(x,y,z)=\frac{1}{2}(\tau_1x^2+\tau_2y^2+\tau_3z^2)+\varepsilon^5\omega(\varepsilon^{-2}(x,y,z)),
\end{equation}
where $\varepsilon$ and $\tau_i$  are arbitrary constants, and $\omega$ is a given smooth function.
\begin{theorem}
Assume that  $f\in\textit{C}^{\infty}(\Omega\times\mathbb{R}\times\mathbb{R}^3)$, then for any $Z_0=(x_0,u_0,p_0)\in\Omega\times\mathbb{R}\times\mathbb{R}^3$, we have that 

 $(1)$ If $f(Z_0)=0$, then  \eqref{2hesi} admits a $1$-convex $\textit{C}^{\infty}$ local solution which is not convex.

   $(2)$ If $f\geq0$ near $Z_0$, then  \eqref{2hesi} admits a $2$-convex $\textit{C}^{\infty}$ local solution which is not convex. 
   
   If $f(Z_0)>0$,  \eqref{2hesi} admits a convex $\textit{C}^{\infty}$ local solution.

   $(3)$ If $f(Z_0)<0$, then  \eqref{2hesi} admits a $1$-convex $\textit{C}^{\infty}$ local solution which is not $2$-convex. \\
Moreover, the equation \eqref{2hesi} is uniformly elliptic with respect to the above local solutions.
\end{theorem}

In this part we obtain the one-parameter groups generated by some operators, and since these groups are symmetry groups of $\mathrm{HESI}$ equation, then the above solution  transforms to another solution of $\mathrm{HESI}$ equation.

\begin{itemize}
\item\textbf{Case 1:} The operator $V_1=\partial_u$  that produces one-parameter group $G_1=(x,y,z, t+u)$, $t\in\mathbb{R}$, for equation ${\cal S}_2[u]=f(x,y,z)$, transforms solution\eqref{solution} to the following solution
\begin{align*}
u(x)=\frac{1}{2}(\tau_1x^2+\tau_2y^2+\tau_3z^2)-t+\varepsilon^5\omega(\varepsilon^{-2}(x,y,z)),
\end{align*}
where $\varepsilon,\tau_i\in\mathbb{R}$.
\item\textbf{Case 2:} The operator $V_2=x\partial_u$  that produces one-parameter group $G_2=(x,y,z, tx+u)$, $t\in\mathbb{R}$, for equation ${\cal S}_2[u]=f(x,y,z)$, transforms solution\eqref{solution} to the following solution
\begin{align*}
u(x)=\frac{1}{2}(\tau_1x^2+\tau_2y^2+\tau_3z^2)-tx+\varepsilon^5\omega(\varepsilon^{-2}(x,y,z)),
\end{align*}
where $\varepsilon,\tau_i\in\mathbb{R}$.
In such a manner are the following cases:
\item\textbf{Case 3:} The operator $V_3=y\partial_u$, $G_3=(x,y,z, ty+u)$, $t\in\mathbb{R}$, for equation ${\cal S}_2[u]=f(x,y,z)$; So
\begin{align*}
u(x)=\frac{1}{2}(\tau_1x^2+\tau_2y^2+\tau_3z^2)-ty+\varepsilon^5\omega(\varepsilon^{-2}(x,y,z)), 
\end{align*}
where $\varepsilon,\tau_i\in\mathbb{R}$.
\item\textbf{Case 4:} The operator $V_4=z\partial_u$, $G_4=(x,y,z, tz+u)$, $t\in\mathbb{R}$, for equation ${\cal S}_2[u]=f(x,y,z)$; So
\begin{align*}
u(x)=\frac{1}{2}(\tau_1x^2+\tau_2y^2+\tau_3z^2)-tz+\varepsilon^5\omega(\varepsilon^{-2}(x,y,z)),
\end{align*}
where $\varepsilon,\tau_i\in\mathbb{R}$.
\item\textbf{Case 5:} The operator $V_5=\pm\partial_x +u\partial_u$, $G_5=(x\pm t,y,z, tu+u)$, $t\in\mathbb{R}$,\\for equation ${\cal S}_2[u]=exp(\pm2x)H(y,z)$; So
\begin{align*}
u(x)&=\frac{1}{2(t+1)}\big[\tau_1(x\pm{t})^2+\tau_2y^2+\tau_3z^2+2\varepsilon^5\omega(\varepsilon^{-2}(x\pm t,y,z))\big],
\end{align*}
where $\varepsilon,\tau_i\in\mathbb{R}$.
\item\textbf{Case 6:} The operator $V_5=\gamma_1z\partial_y -\gamma_1y\partial_z +u\partial_u$, \\$G_5=(x,z\sin(\gamma_1t)+y\cos(\gamma_1t),z\cos(\gamma_1t)-y\sin(\gamma_1t), tu+u)$, $t\in\mathbb{R}$,\\
 for equation ${\cal S}_2[u]=\di\exp(\frac{2}{\gamma_1}{\tan}^{-1}(\frac{y}{z})) H(x,y^2+z^2)$; So
\begin{align*}
u(x)&=\frac{1}{2(t+1)}\big[\tau_1x^2+\tau_2(z\sin(\gamma_1t)+y\cos(\gamma_1t))^2+\tau_3(z\cos(\gamma_1t)-y\sin(\gamma_1t))^2\\
&+2\varepsilon^5\omega(\varepsilon^{-2}(x,z\sin(\gamma_1t)+y\cos(\gamma_1t),z\cos(\gamma_1t)-y\sin(\gamma_1t)))\big],
\end{align*}
where $\varepsilon,\tau_i\in\mathbb{R}$.
\item\textbf{Case 7:} The operator $V_5=\pm\partial_x +\gamma_2z\partial_y -\gamma_2y\partial_z +u\partial_u$,\\
$G_5=(x\pm t,z\sin(\gamma_2t)+y\cos(\gamma_2t),z\cos(\gamma_2t)-y\sin(\gamma_2t), tu+u)$, $t\in\mathbb{R}$,\\ for equation ${\cal S}_2[u]=\di \exp(\frac{2}{\gamma_2}{\tan}^{-1}(\frac{y}{z}))\di H(y^2+z^2,x\mp \frac{1}{\gamma_2}{\tan}^{-1}(\frac{y}{z}))$; So
\begin{align*}
u(x)&=\frac{1}{2(t+1)}\big[\tau_1(x\pm t)^2+\tau_2(z\sin(\gamma_2t)+y\cos(\gamma_2t))^2+\tau_3(z\cos(\gamma_2t)-y\sin(\gamma_2t))^2\\
&+2\varepsilon^5\omega(\varepsilon^{-2}(x\pm t,z\sin(\gamma_2t)+y\cos(\gamma_2t),z\cos(\gamma_2t)-y\sin(\gamma_2t)))\big],
\end{align*}
where $\varepsilon,\tau_i\in\mathbb{R}$.
\item\textbf{Case 8:} The operator $V_5=\alpha_1z\partial_x -\alpha_1x\partial_z +u\partial_u $,\\ $G_5=(z\sin(\alpha_1t)+x\cos(\alpha_1t),y,z\cos(\alpha_1t)-x\sin(\alpha_1t), tu+u)$, $t\in\mathbb{R}$,\\
 for equation ${\cal S}_2[u]=\exp(\frac{2}{\alpha_1}{\tan}^{-1}(\frac{x}{z}))H(y,x^2+z^2)$; So
\begin{align*}
u(x)&=\frac{1}{2(t+1)}\big[\tau_1(z\sin(\alpha_1t)+x\cos(\alpha_1t))^2+\tau_2y^2+\tau_3(z\cos(\alpha_1t)-x\sin(\alpha_1t))^2\\
&+2\varepsilon^5\omega(\varepsilon^{-2}(z\sin(\alpha_1t)+x\cos(\alpha_1t),y,z\cos(\alpha_1t)-x\sin(\alpha_1t)))\big],
\end{align*}
where $\varepsilon,\tau_i\in\mathbb{R}$.
\item\textbf{Case 9:} The operator $V_5=\pm\partial_y +u\partial_u$, $G_5=(x,y\pm t,z, tu+u)$, $t\in\mathbb{R}$,\\ for equation ${\cal S}_2[u]=exp(\pm2y)H(x,z)$; So
\begin{align*}
u(x)=\frac{1}{2(t+1)}\big[\tau_1x^2+\tau_2(y\pm t)^2+\tau_3z^2+2\varepsilon^5\omega(\varepsilon^{-2}(x,y\pm t,z))\big],
\end{align*}
where $\varepsilon,\tau_i\in\mathbb{R}$.
\item\textbf{Case 10:} The operator $V_5= \pm\partial_y +\alpha_2z\partial_x -\alpha_2x\partial_z +u\partial_u$,\\ $G_5=(z\sin(\alpha_2t)+x\cos(\alpha_2t),y\pm t,z\cos(\alpha_2t)-x\sin(\alpha_2t), tu+u)$, $t\in\mathbb{R}$,\\
 for equation ${\cal S}_2[u]=\di \exp(\frac{2}{\alpha_2}{\tan}^{-1}(\frac{x}{z}))\di H(x^2+z^2,y\mp \frac{1}{\alpha_2}{\tan}^{-1}(\frac{x}{z}))$; So
\begin{align*}
u(x)&=\frac{1}{2(t+1)}\big[\tau_1(z\sin(\alpha_2t)+x\cos(\alpha_2t))^2+\tau_2(y\pm t)^2+\tau_3(z\cos(\alpha_2t)-x\sin(\alpha_2t))^2\\
&+2\varepsilon^5\omega(\varepsilon^{-2}(z\sin(\alpha_2t)+x\cos(\alpha_2t),y\pm t,z\cos(\alpha_2t)-x\sin(\alpha_2t)))\big],
\end{align*}
where $\varepsilon,\tau_i\in\mathbb{R}$.
\item\textbf{Case 11:} The operator $V_5=\beta_1y\partial_x -\beta_1x\partial_y +u\partial_u $,\\
$G_5=(y\sin(\beta_1t)+x\cos(\beta_1t),y\cos(\beta_1t)-x\sin(\beta_1t),z, tu+u)$, $t\in\mathbb{R}$,\\
 for equation ${\cal S}_2[u]= \exp(\frac{2}{\beta_1}{\tan}^{-1}(\frac{x}{y}))H(z,x^2+y^2)$; So
\begin{align*}
u(x)&=\frac{1}{2(t+1)}\Big[\tau_1(y\sin(\beta_1t)+x\cos(\beta_1t))^2+\tau_2(y\cos(\beta_1t)-x\sin(\beta_1t))^2+\tau_3z^2\\
&+2\varepsilon^5\omega(\varepsilon^{-2}(y\sin(\beta_1t)+x\cos(\beta_1t),y\cos(\beta_1t)-x\sin(\beta_1t),z))\Big],
\end{align*}
where $\varepsilon,\tau_i\in\mathbb{R}$.
\item\textbf{Case 12:} The operator $V_5=\pm\partial_z +u\partial_u$,  $G_5=(x,y,z\pm t, tu+u)$, $t\in\mathbb{R}$,\\ for equation ${\cal S}_2[u]= {\exp(\pm 2z)}{H(x,y)}$; So
\begin{align*}
u(x)=\frac{1}{2(t+1)}\Big[\tau_1x^2+\tau_2y^2+\tau_3(z\pm t)^2+2\varepsilon^5\omega(\varepsilon^{-2}(x,y,z\pm t))\Big],
\end{align*}
where $\varepsilon,\tau_i\in\mathbb{R}$.
\item\textbf{Case 13:} The operator $V_5=\pm\partial_z +\beta_2y\partial_x -\beta_2x\partial_y +u\partial_u$,\\  $G_5=(y\sin(\beta_2t)+x\cos(\beta_2t),y\cos(\beta_2t)-x\sin(\beta_2t),z\pm t, tu+u)$, $t\in\mathbb{R}$,\\
 for equation ${\cal S}_2[u]=\di \exp(\frac{2}{\beta_2}{\tan}^{-1}(\frac{x}{y}))\di H(x^2+y^2,z\mp \frac{1}{\beta_2}{\tan}^{-1}(\frac{x}{y}))$; So
\begin{align*}
u(x)&=\frac{1}{2(t+1)}\big[\tau_1(y\sin(\beta_2t)+x\cos(\beta_2t))^2+\tau_2(y\cos(\beta_2t)-x\sin(\beta_2t))^2+\tau_3(z\pm t)^2\\
&+2\varepsilon^5\omega(\varepsilon^{-2}(y\sin(\beta_2t)+x\cos(\beta_2t),y\cos(\beta_2t)-x\sin(\beta_2t),z\pm t))\big],
\end{align*}
where $\varepsilon,\tau_i\in\mathbb{R}$.
\item\textbf{Case 14:} The operator $V_5=x\partial_x+y\partial_y+z\partial_z$,\\
$G_5=(e^tx,e^ty,e^tz,u)$, $t\in\mathbb{R}$, for equation ${\cal S}_2[u]= x^{2c-4}H(\frac{y}{x},\frac{z}{x})$; So
\begin{align*}
u(x)=\frac{1}{2}(\tau_1e^{2t}x^2+\tau_2e^{2t}y^2+\tau_3e^{2t}z^2)+\varepsilon^5\omega(\varepsilon^{-2}(e^tx,e^ty,e^tz)),
\end{align*}
where $\varepsilon,\tau_i\in\mathbb{R}$.
\item\textbf{Case 15:} The operator $V_5=x\partial_x+(y\pm 1)\partial_y+z\partial_z$, \\
$G_5=(e^tx,e^ty\pm e^t\mp1,e^tz,u)$, $t\in\mathbb{R}$, for equation ${\cal S}_2[u]= x^{2c-4}H(\frac{y\pm 1}{x},\frac{z}{x})$; So
\begin{align*}
u(x)&=\frac{1}{2}(\tau_1e^{2t}x^2+\tau_2(e^ty\pm e^t\mp1)^2+\tau_3e^{2t}z^2)+\varepsilon^5\omega(\varepsilon^{-2}(e^tx,e^ty\pm e^t\mp1,e^tz)),
\end{align*}
where $\varepsilon,\tau_i\in\mathbb{R}$.
\end{itemize}
\section{Conclusion}
In this paper,  we performed preliminary group classification on equation, by studying the class of $3$-dimensional nonlinear $2-$hessian equations ${\cal S}_2[u]=f(x,y,z)$, and investigating the algebraic structure of the symmetry groups for the equation. Then, we obtained an optimal system of one-dimensional Lie subalgebras of this equation, with the aid of propositions \eqref{proposition1} and \eqref{proposition2}. The result of these work is a wide class of equations which summarized in table \ref{tab3}. And at the end of this work, some exact solutions of $2-$hessian equation are presented.
 Of course, the results in table \ref{tab3} can be continued for remainder vectors, and it is possible to obtain the corresponding reduced equations for all the cases in the classification in up comming works.

\end{document}